\def\1{1\!\!1} 
\documentclass{conm-p-l}

\newtheorem{theorem}{Theorem}[section]
\newtheorem{lemma}[theorem]{Lemma}

\theoremstyle{definition}

\theoremstyle{remark}
\newtheorem{remark}[theorem]{Remark}

\numberwithin{equation}{section}

\def\N{I\!\!N}                      \def\R{I\!\!R}

\def\1{1\!\!1}  

\def\and{\text{ and }}

\def\h{{\text h}}            
\def\hmu{\h_\mu}           

\def\H{{\text H}}             
                                        
       \def\ep{\text{e}}
         \def\P{\text{{\rm P}}}

\def\a{\alpha}                \def\b{\beta}              
                         \def\f{\phi}
\def\g{\gamma}                           \def\l{\lambda} 
              \def\om{\omega}           
\def\Sg{\Sigma}               \def\sg{\sigma}

\def\bi{\bigcap}              \def\bu{\bigcup}          
\def\({\bigl(}                \def\){\bigr)}           
\def\lt{\left}                \def\rt{\right}

\def\ld{\ldots}                        \def\^{\tilde}

\def\es{\emptyset}            \def\sms{\setminus}
\def\sbt{\subset}

\def\comp{\asymp}
           
      \def\fr{\noindent}

\def\om{\omega}

\begin{document}

\title{A new class of positive recurrent functions}

\author{Pawel Hanus}
\address{Department of Mathematics, University of North Texas, Denton, 
TX 76203-5118, USA}
\email{pgh0001@jove.acs.unt.edu, mauldin\@unt.edu}
\author{Mariusz Urba\'nski}
\address{Department of Mathematics, University of North Texas, Denton, 
TX 76203-5118, USA}
\email{urbanski@unt.edu}
\thanks{The second author was supported in part by NSF Grant DMS 9801583}

\date{June 18, 1999}
\begin{abstract} In [Sa] Sarig has introduced and explored the
concept of positively recurrent functions. In this paper we construct
a natural wide class of such functions and we show
that they have stronger ergodic properties than the general functions 
considered in [Sa].
\end{abstract}

\maketitle

\section{Preliminaries} 

In [Sa] Sarig has introduced and explored the
concept of positively recurrent functions. In this paper, using the
concept of an iterated function system,  we construct
a natural wide class of positively recurrent functions and we show
that they have stronger properties than the general functions considered
in [Sa]. In some parts our exposition is similar and follows the
approach developed in [MU1] and [Wa], where also the idea of embedding 
the infinite dimensional shift space into a compact metric space and the 
Shauder-Tichonov fixed-point theorem have been used. 
To begin with, let $\N$ be the set of
positive integers
and let $\Sg=\N^\infty$ be the infinitely dimensional shift space
equipped with the product topology. Let $\sg:\Sg\to\Sg$ be the shift
transformation (cutting out the first coordinate),
$\sg(\{x_n\}_{n=1}^\infty) =(\{x_n\}_{n=2}^\infty)$. Fix $\b>0$. If
$\phi:\Sg\to \R$ and $n\ge 1$, we set
$$
V_n(\phi)=\sup\{|\phi(x)-\phi(y)|:x_1=y_1,x_2=y_2,\ld, x_n=y_n\}.
$$
The function $\phi$ is said to be H\"older continuous of order $\b$ if
and only if 
$$
V(\phi)=\sup_{n\ge 1}\{\ep^{\b n}V_n(\phi)\}<\infty.
$$
We also assume that
\begin{equation}
\sup_{\om\in\Sg}\sum_{\tau\in\sg^{-1}(\om)}\ep^{\phi(\tau)} <\infty.
\end{equation}
This assumption allows us to introduce the Perron-Frobenius-Ruelle
operator ${\it L}_\phi:C_b(\Sg)\to C_b(\Sg)$,
$$
{\it L}_\phi(g)(\om)=\sum_{\tau\in\sg^{-1}(\om)}\ep^{\phi(\tau)}g(\tau)
$$
acting on $C_b(\Sg)$, the space of all bounded continuous real-valued
functions on $\Sg$ equipped with the norm $||\cdot||_0$, where
$||k||_0=\sup_{x\in\Sg}|k(x)|$. Moreover,
$$
||{\it L}_\phi||_0\le{\it L}_\phi(\1)=\sup_{\om\in\Sg}
\sum_{\tau\in\sg^{-1}(\om)}\ep^{\phi(\tau)}<\infty.
$$
We extend the standard definition of topological pressure by setting
\begin{equation}
\P(\phi)=\lim_{n\to\infty}{1\over n}\log\lt(\sum_{|\om|=n}\sup_{\tau\in
[\om]}\exp\lt(\sum_{j=0}^{n-1}\phi\circ\sg^j(\tau)\rt)\rt), 
\end{equation}
where $[\om]=\{\rho\in\Sg:\rho_1=\om_1,\rho_2=\om_2,\ld,\rho_{|\om|}=
\om_{|\om|}\}$. Notice that the limit exists since the partition functions
$$
Z_n(\phi)
=\sum_{|\om|=n}\sup_{\tau\in [\om]} \exp\lt(\sum_{j=0}^{n-1}\phi
\circ\sg^j(\tau)\rt)
$$
form a subadditive sequence. Notice also that our
definition of pressure formally differs from that provided by Sarig in
[Sa] which reads that given $i\in \N$
\begin{equation}
\P(\phi)=\lim_{n\to\infty}{1\over n}\log Z_n(\phi,i), 
\end{equation}
where
$$
Z_n(\phi,i)=\sum\exp\lt(\sum_{j=0}^{n-1}\phi\circ\sg^j(\om)\rt)
$$
and the summation is taken over all elements $\om$ satisfying $\sg^n(\om)
=\om$ and $\om_1=i$. However in [Sa] Sarig proves Theorem 2 which says
that $\P(\phi)=\sup\{\P(\phi|_Y)\}$, where the supremum is taken over
all topologically mixing subshifts of finite type $Y\sbt \Sg$ and the
same proof goes through with (1.3) replaced by (1.2) (comp. Theorem 3.1 
of [MU2]). Thus we have the following.

\begin{lemma} The definitions of topological pressures given by
(1.2) and (1.3) coincide.
\end{lemma}
\

\fr Here is a direct proof of this lemma communicated to us by Sarig
Omri: Fix $i\in \N$. Using H\"older continuity of the function $\phi$ we
can write 
$$
Z_n(\phi)
\comp \sum_{|\om|=n}\exp\(\sum_{j=0}^{n-1}\phi(\om^\infty)\)
\comp \sum_{|\om|=n}\exp\(\sum_{j=0}^n\phi((i\om)^\infty)\)
=Z_{n+1}(\phi,i).
$$
Thus the lemma is proved. 

\

\fr Following the definition 2 of [Sa] we call the function $\phi:\Sg\to\R$
positive recurrent if for every $i\in \N$ there exists a constant $M_i$
and an integer $N_i$ such that for all $n\ge N_i$
$$
M_a^{-1}\le Z_n(\phi,i)\l^{-n}\le M_i
$$
for some $\l>0$. As we already have said the main purpose of this paper
is to provide a wide natural class of examples of positive recurrent
potential which additionally satisfy much stronger properties than those
claimed in Theorem 4 of [Sa]. In order to describe our setting let
$(X,d)$ be a compact metric space and let $\phi_i:X\to X$, $i\in \N$, be
a family of uniform contractions, i.e. $d(\phi_i(x),\phi_i(y))\le sd(x,y)$ for
all $i\in \N$, $x,y\in X$ and some $s<1$. Given $\om\in\Sg$ consider the
intersection $\bi_{n\ge 1}\phi_{\om|_n}(X)$, where $\phi_{\om|_n}=
\phi_{\om_1}\circ \ld\circ \phi_{\om_n}$. Since $\phi_{\om|_n}(X)$,
$n\ge 1$, form a descending family of compact sets, this intersection is
non-empty and since the maps $\phi_i$, $i\in \N$, are uniform
contractions, it is a singleton. So, we have defined a projection map
$\pi:\Sg\to X$ given by the formula
$$
\{\pi(\om)\}=\bi_{n\ge 1}\phi_{\om|_n}(X).
$$
$J$, the range of $\pi$, is said to be the limit set of the iterated
function system $\phi_i:X\to X$, $i\in \N$. Let now $\phi^{(i)}:X\to 
\R$, $i\in\N$, be a family of continuous functions such that
\begin{equation}
\sup_X\sum_{i\in \N}\ep^{\phi^{(i)}(x)}<\infty. 
\end{equation}
We define a function $\phi:\Sg\to \R$ by setting
\begin{equation}
\phi(\om)=\phi^{(\om_1)}(\pi(\sg(\om))). 
\end{equation}
It easily follows from (1.4) that $\P(\phi)<\infty$. In the next section we
shall prove the following.

\begin{theorem} Suppose that the function $\phi:\Sg\to \R$ defined by (1.5) 
and satisfying (1.4) is H\"older continuous. Let ${\it L}_\phi^*$ be the
operator conjugate to ${\it L}_\phi$. Then $\phi$ is positive recurrent with
$\l=\ep^{\P(\phi)}$. Moreover there exists $M>0$ such that $M^{-1}\le
\l^{-n}{\it L}_\phi^n(\1)\le M$ for all $n\ge 1$. Suppose  additionally
that $\phi_i(X)\cap\phi_j(X)=\es$ for all $i,j\in \N$, $i\ne j$. Then
there are a probability measure $\nu$ on $\Sg$ and a bounded away from
zero and infinity H\"older continuous function $h:\Sg\to (0,\infty)$
such that ${\it L}_\phi^*(\nu)=\l\nu$, ${\it L}_\phi(h)=\l h$, $\nu(h)=1$ and
$\l^{-n} {\it L}_\phi^n(g)\to (\int gd\nu)h$ uniformly for every
uniformly continuous bounded function $g$. Additionally 
$\l^{-n}{\it L}_\phi^n(g)
\to (\int gd\nu)h$ exponentially fast for each H\"older continuous bounded 
function $g$.
\end{theorem}

\section{Proof of Theorem~1.2.} Define first an auxiliary
Perron-Frobenius operator $\^L_\phi:C(X)\to C(X)$ given by the formula
$$
\^L_\phi(g)(x)=\sum_{i\in\N}\ep^{\phi^{(i)}(x)}g(\phi_i(x)).
$$
$\^L_\phi$ is continuous, positive and $||\^L_\phi||_0\le
\sup_X\sum_{i\in\N}\ep^{\phi^{(i)}(x)}<\infty$. Let $\^L_\phi^*:C(X)^*\to
C(X)^*$ be the conjugate operator and following Bowen's approach from [Bo]
consider the map
$$
\mu\mapsto {\^L_\phi^*(\mu)\over \^L_\phi^*(\mu)(\1)}.
$$
of the space of Borel probability measures on $X$ into itself. This map
is continuous in the weak-* topology of measures and therefore, in view
of the Schauder-Tichonov theorem, it has a fixed point, say $m_\phi$. Thus
\begin{equation}
\^L_\phi^*(m_\phi)=\l m_\phi 
\end{equation}
with $\l=\^L_\phi^*(m_\phi)(\1)$.

\

\fr Given $n\ge 1$ and $\om\in \N^n$, denote $\sum_{j=1}^n\phi^{(\om_j)}
\circ \phi_{\sg^j\om}$ by $S_\om(\phi)$. Let us then prove the following.

\begin{lemma} If $x,y\in \phi_\tau(X)$ for some $\tau\in I^*$, then
for all $\om\in I^*$ 
$$
|S_\om(\phi)(x)-S_\om(\phi)(y)|\le
{V(\phi)\over 1-\ep^{-\b}}\ep^{-\b|\tau|} 
$$
\end{lemma}
\begin{proof} Let $n=|\om|$. Write $x=\phi_\tau(u)$, $y=\phi_\tau(w)$,
where $u,w\in X$. By (2.1) we get 
$$\aligned
\big|\sum_{j=1}^n\phi^{(\om_j)}(\phi_{\sg^j\om}(x)) &-
\sum_{j=1}^n\phi^{(\om_j)}(\phi_{\sg^j\om}(y))\big|= \\
&=\left|\sum_{j=1}^n\phi^{(\om\tau)_j}\circ\phi_{\sg^j\om\tau}(u)
-\sum_{j=1}^n\phi^{(\om\tau)_j}\circ \phi_{\sg^j\om\tau}(w)\right| \\
&\le \sum_{j=1}^n\left|\phi^{(\om\tau)_j}\circ\phi_{\sg^j\om\tau}(u) - 
\phi^{(\om\tau)_j}\circ\phi_{\sg^j\om\tau}(w)\right| \\
&\le \sum_{j=1}^nV(\phi)\ep^{-\b(n+|\tau|-j)} \\
&\le {V(\phi)\over 1-\ep^{-\b}}\ep^{-\b |\tau|} \endaligned
$$
The proof is finished. 
\end{proof}

\begin{remark} 
We allow in Lemma~2.1 $\tau$ to be the empty word
$\es$. Then $\phi_\es=\text{Id}_X$ and $|\es|=0$.
\end{remark}

\fr Set
$$
Q=\exp\lt(V(\phi){\ep^{-\b}\over 1-\ep^{-\b}}\rt).
$$
We shall prove the following.

\begin{lemma} The eigenvalue $\l$ (see 2.1) of the dual
Perron-Frobenius operator is equal to $\ep^{\P(\phi)}$.
\end{lemma}
\begin{proof} 
Iterating (2.1) we get 
$$\aligned
\l^n
&=\l^nm_\phi(\1)=\^L_\phi^{*n}(\1)=\int_X\^L_\phi^n(\1)dm_\phi \\
&=\int_X\sum_{|\om|=n}\exp(S_\om(\phi)(x))
\le \sum_{|\om|=n}||\exp(S_\om(\phi))||_0. \endaligned
$$
So,
$$
\log\l
\le \lim_{n\to\infty}{1\over n}\log \lt(\sum_{|\om|=n}||
\exp(S_\om(\phi))||_0\rt) =\P(\phi).
$$
Fix now $\om\in I^n$ and take a point $x_\om$ where the function
$S_\om(\phi)$ takes on its maximum. In view of Lemma~2.1, for every
$x\in X$ we have
$$
\sum_{|\om|=n}\exp(S_\om(\phi)(x))
\ge Q^{-1} \sum_{|\om|=n}\exp(S_\om(\phi)(x_\om))
=Q^{-1} \sum_{|\om|=n}||\exp(S_\om(\phi))||_0.
$$
Hence, iterating (2.1) as before,
$$
\l^n
=\int_X\sum_{|\om|=n}\exp(S_\om(\phi))dm_\phi
\ge Q^{-1}\sum_{|\om|=n}||\exp(S_\om(\phi))||_0.
$$
So, $\log\l\ge \lim_{n\to\infty}{1\over n}\log\(
\sum_{|\om|=n}||\exp(S_\om(\phi))||_0\) =\P(\phi)$. The proof is
finished. 
\end{proof}

\fr Let $\^L_0$ and ${\it L}_0$ denote the corresponding normalized 
Perron-Frobenius operators, i.e. $\^L_0= \ep^{-\P(\phi)}\^L_\phi$ and
${\it L}_0=\ep^{-\P(\phi)}{\it L}_\phi$. We shall prove the following.

\begin{theorem} 
$m_\phi(J)=1$.
\end{theorem}
\begin{proof} 
Since by (2.1)
\begin{equation}
\^L_0^*(m_\phi)=m_\phi 
\end{equation}
and consequently $\^L_0^{*n}(m_\phi)=m_\phi$ for all $n\ge 0$, we have 
\begin{equation}
\int_X\sum_{|\om|=n}\exp\(S_\om(\phi)-\P(\phi)n\)\cdot (f\circ
\phi_\om)dm_\phi =\int_Xfdm_\phi 
\end{equation}
for all $n\ge 0$ and all continuous functions $f:X\to\R$. Since this
equality extends to all bounded measurable functions $f$, we get
\begin{equation}
\aligned
m_\phi(A)
&=\sum_{\tau\in I^n}\int
\exp\(S_\tau(\phi)-\P(\phi)n\)\cdot\1_{\phi_\om(A)}\circ \phi_\tau dm_\phi \\
&\ge \int_A\exp\(S_\om(\phi)-\P(\phi)n\)dm_\phi 
\endaligned
\end{equation}
for all $n\ge 0$, all $\om\in I^n$, and all Borel sets $A\sbt X$. Now,
for each $n\ge 1$ set $X_n=\bu_{|\om|=n}\phi_\om(X)$. Then $\1_{X_n}
\circ\phi_\om=\1$ for all $\om\in \N^n$. Thus apllying (2.3) to the
function $f=\1_{X_n}$ and later to the function $f=\1$, we obtain
$$\aligned
m_\phi(X_n)
&=\int_X\sum_{|\om|=n}\exp\(S_\om(\phi)-\P(\phi)n\)\cdot (\1_{X_n}\circ
\phi_\om)dm_\phi \\
&= \int_X\sum_{|\om|=n}\exp\(S_\om(\phi)-\P(\phi)n\)dm_\phi
= \int \1 dm_\phi 
=1. \endaligned 
$$
Hence $m_\phi(J)=m_\phi\(\bi_{n\ge 1}X_n\)=1$. The proof is
complete. 
\end{proof}

\begin{theorem} For all $n\ge 1$
$$
Q^{-1} \le \^L_0^n(\1) \le Q.
$$
\end{theorem}
\begin{proof} 
Given $n\ge 1$ by (2.3) there exits $x_n\in X$ such that 
$\^L_0^n(\1)(x_n)\le 1$. It then follows from Lemma~2.1 that for every
$x\in X$, $\^L_0^n(\1)\le Q$. Similarly by
(2.3) there exists $y_n\in X$ such that $\^L_0^n(\1)\ge 1$. It then
follows from Lemma~2.1  that for every $x\in X$, $\^L_0^n(\1)
\ge Q^{-1}$. The proof is finished. 
\end{proof}

\

\fr So far we have worked downstairs in the compact space $X$. It is now
time to lift our considerations up to the shift space $\Sg$.

\begin{lemma} There exists a unique Borel probability measure
$\^m_\phi$ on $\N^\infty$ such that $\^m_\phi([\om])=\int \exp
\(S_\om(\phi)-\P(\phi)n\)dm_\phi$ for all $\om\in \N^*$.
\end{lemma}
\begin{proof} In view of (2.4)
$\int\exp\(S_\om(\phi)-\P(\phi)n\)dm_\phi =1$ for all $n\ge 1$ and
therefore one can define a Borel probability measure $m_n$ on $C_n$,
the algebra generated by the cylinder sets of the form $[\om]$,
$\om\in \N^n$, putting $m_n([\om])=\int\exp\(S_\om(\phi) 
-\P(\phi)n\)dm_\phi$. Hence, applying (2.4) again we get for all $\om\in
\N^n$.
$$\aligned
m_{n+1}(\om)
&=\sum_{i\in \N}m_{n+1}([\om i])
 =\sum_{i\in \N} \int\exp\(S_{\om i}(\phi)-\P(\phi)n\)dm_\phi \\
&=\int \sum_{i\in \N} \exp\lt(\sum_{j=1}^n\phi^{(\om_j)}\circ
\phi_{\sg^j(\om i)}-\P(\phi)n + \phi^{(i)}-\P(\phi)\rt)dm_\phi \\
&= \int \sum_{i\in \N} \exp\(S_\om\circ \phi_i -\P(\phi)n\) 
\exp\(\phi^{(i)}-\P(\phi)\)dm_\phi \\
&=\int \^L_0\(\exp(S_\om(\phi)-\P(\phi))\)dm_\phi
 = \int\exp\(S_\om(\phi)-\P(\phi)\)dm_\phi
 = m_n([\om]) \endaligned
$$
and therefore in view of Kolmogorov's extension theorem there exists a
unique probability measure $\^m_\phi$ on $\N^\infty$ such that 
$\^m_\phi([\om])=
\^m_{|\om|}([\om])$ for all $\om \in \N^*$. The proof is complete. \end{proof}

\

\fr Now we are ready to prove that the function $\phi$ is positive
recurrent. Let us first notice that
$$\aligned
{\it L}_\phi(\1)(\om)
&=\sum_{\tau\in\sg^{-1}(\om)}\ep^{\phi(\tau)}
 =\sum_{\tau\in\sg^{-1}(\om)}\exp\(\phi^{(\tau_1)}(\pi(\sg(\tau)))\) \\
&=\sum_{\tau\in\sg^{-1}(\om)}\exp\(\phi^{(\tau_1)}(\pi(\om))\)
 =\sum_{i\in\N}\ep^{\phi^{(i)}(\pi(\om))}
 =\^L_\phi(\1)(\pi(\om)). \endaligned
 $$
Since $\^L_0=\ep^{-\P(\phi)}\^L_\phi$, it then follows from Theorem~2.4
that as $M$ we can take $Q$. In order to demonstrate that the function
$\phi$ is positive recurrent we first show that
$$
{Z_n(\phi,i)\over {\it L}_\phi^n(\1)(\om)} \le M_i
$$
for all $n\ge 1$, $\om\in\Sg$, and some constant $M_i>0$. So fix 
$\om\in\Sg$. We shall define an injection $j$ from $\{\rho\in\Sg:
\sg^n(\rho)=\rho \text{ and } \rho_1=i\}$ into $\sg^{-n}(\om)$ as
follows: $j(\rho)=\rho_1\rho_2\ld\rho_n\om$. Now, by Lemma~2.1
$$
\left|\sum_{j=0}^{n-1}\phi(\sg^j(\rho))
-\sum_{j=0}^{n-1}\phi(\sg^j(j(\rho)))\right| \le \log Q
$$
and therefore $Z_n(\phi,i)\le Q{\it L}_\phi^n(\1)(\om)$. Thus by Theorem~2.4
and the definition of the operators $\^L_0$ and ${\it L}_0$, $Z_n(\phi,i)\le
M_i\l^n$, where $M_i=Q^2$. Now we shall prove that  $Z_n(\phi,i)\ge
M_a'\l^n$ for some constant $M_i'$ and all $n\ge 1$. We demonstrate
first that for all $n\ge 1$ and all $i\in \Sg$
$$
{\it L}_0(\1_{[i]}) \ge \^m_\phi([i]).
$$
Indeed, since $\int {\it L}_0(\1_{[i]})d\^m_\phi= \int \1_{[i]})d\^m_\phi= 
\^m_\phi([i]) >0$, there exists $\tau\in\Sg$ such that ${\it L}_0(\1_{[i]})
(\tau)\ge \^m_\phi([i])$. It the follows from Lemma~2.1 that for every
$\om\in\Sg$
$$\aligned
{\it L}_0^n(\1_{[i]})(\om)
&=\sum_{\rho\in\sg^{-n}(\om)}\exp\lt(\sum_{j=0}^{n-1}\phi\circ\sg^j(\rho)
  \1_{[i]}(\rho)\rt) \\
&\ge Q^{-1}\sum_{\rho\in\sg^{-n}(\tau)}\exp\lt(\sum_{j=0}^{n-1}\phi\circ
\sg^j(\rho)  \1_{[i]}(\rho)\rt)
 =Q^{-1} {\it L}_0(\1_{[i]})(\tau) \\
&\ge  \^m_\phi([i]). \endaligned 
$$
Hence ${\it L}_\phi^n(\1_{[i]})(\om)\ge \l^n \^m_\phi([i])$. So, in order to
conclude the proof that $\phi$ is positively recurrent it suffices now
to show that
$$
{Z_n(\phi,i)\over {\it L}_\phi^n(\1_{[i]})(\om)} \ge M_i''
$$
for all $n\ge 1$, all $\om\in\Sg$ and some constant $M_i''>0$. Indeed, we
shall define an injection $k$ from $\sg^{-n}(\om)\cap [i]$ to $\{\rho:
\Sg: \sg^n(\rho)=\rho \text{ and } \rho_1=i\}$ by taking as $k(\tau)$
the infinite concatenation of the first $n$ words of $\tau$. Then by
Lemma~2.1,
$$
\left|\sum_{j=0}^{n-1}\phi(\sg^j(\tau))
-\sum_{j=0}^{n-1}\phi(\sg^j(k(\tau)))\right| \le \log Q
$$
and therefore
$$\aligned
{\it L}_\phi^n(\1_{[i]})(\om)
&=\sum_{\rho\in\sg^{-n}(\om)}\exp\lt(\sum_{j=0}^{n-1}\phi\circ\sg^j(\rho)
  \1_{[i]}(\rho)\rt) \\
&=\sum_{\rho\in\sg^{-n}(\om)\cap [i]}\exp\lt(\sum_{j=0}^{n-1}\phi\circ
  \sg^j(\rho)\rt) \\
&\le \sum_{\rho\in\sg^{-n}(\om)\cap [i]}\exp\lt(\sum_{j=0}^{n-1}\phi\circ
  \sg^j(k(\rho))+\log Q\rt) \\
&\le Q\sum\exp\lt(\sum_{j=0}^{n-1}\phi
  \circ\sg^j(\rho)\rt) 
  =QZ_n(\phi,i), \endaligned
$$  
where the last summation is taken over all elements $\om$ satisfying
$\sg^n(\om)=\om$ and $\om_1=i$. So, the proof of the positive recurrence
of $\phi$ is complete taking $Q^{-1}$ as $M_i''$. Now we pass to proving
the existence of the measure $\nu$ and the function $h$. We begin with
the following two facts.

\begin{lemma} The measures $m_\phi$ and $\^m_\phi\circ \pi^{-1}$ are
equal. 
\end{lemma}
\begin{proof} Let $A\sbt J$ be an arbitrary closed subset of $J$ and for 
every $n\ge 1$ let $A_n=\{\om\in \N^n:\f_\om(X)\cap A\ne\es\}$. In view of 
(2.3) applied to the characteristic function $\1_A$ we have for all $n\ge 1$
$$\aligned
m_\phi(A) 
&=\sum_{\om\in \N^n}\int\exp\(S_\om(\phi)-\P(\phi)|\om|\)(\1_A\circ\f_\om)
\,dm_\phi \\
&=\sum_{\om\in
   A_n}\int\exp\(S_\om(\phi)-\P(\phi)|\om|\)(\1_A\circ\f_\om)\,dm_\phi  \\
&\le \sum_{\om\in A_n}\int\exp\(S_\om(\phi)-\P(\phi)|\om|\)\,dm_\phi
=\sum_{\om\in A_n}\^m_\phi([\om])
=\^m_\phi\(\bu_{\om\in A_n}[\om]\) \endaligned
$$
Since the family of sets $\{\bu_{\om\in A_n}[\om]:n\ge 1\}$ is descending
and $\bi_{n\ge 1}\bu_{\om\in A_n}[\om]=\pi^{-1}(A)$ we therefore get 
$m_\phi(A)\le\lim_{n\to\infty}\^m_\phi\(\bu_{\om\in A_n}[\om]\)=
\^m_\phi(\pi^{-1}(A))$. Since the limit set $J$ is a metric space ,
using the Baire classification of Borel sets we easily
see that this inequality extends to the family of all Borel subsets of $J$. 
Since both measures $m_\phi$ and $\^m_\phi\circ\pi^{-1}$ are probabilistic
we get $m_\phi=\^m_\phi\circ\pi^{-1}$. The proof is finished. \end{proof}

\

\fr We recall that an invariant mesure of a metric dynamical system is
said to be totally ergodic if it is ergodic with respect to all the
iterates of the system under consideration. 

\begin{theorem} There exists a unique totally ergodic $\sg$-invariant 
probability 
measure $\^\mu_\phi$ absolutely continuous with respect to $\^m_\phi$.
Moreover $\^\mu_\phi$ 
is equivalent with $\^m_\phi$ and $Q^{-1}\le d\^\mu_\phi/d\^m_\phi\le Q$.
\end{theorem}
\begin{proof} First notice 
that, using Lemma~2.5, for each $\om\in \N^*$ and each $n\ge 0$ we have 
$$\aligned
\^m_\phi(\sg^{-n}([\om])) 
&= \sum_{\tau\in \N^n}\^m_\phi([\tau\om])
= \sum_{\tau\in \N^n}
\int\exp\(S_{\tau\om}(\phi)-\P(\phi)|\tau\om|\)\,dm \\
&\ge \sum_{\tau\in \N^n}Q^{-
1}||\exp\(S_\tau(\phi)-\P(\phi)|\tau|\)||_0\exp\(S_\om(\phi-
\P(\phi)|\om|\)\,dm \\ 
&=Q^{-1}\int\exp\(S_\om(\phi-\P(\phi)|\om|\)\,dm_\phi
\sum_{\tau\in \N^n}||\exp\(S_\tau(\phi-\P(\phi)|\tau|\)||_0 \\
&\ge Q^{-1}\^m_\phi([\om])\^m_\phi(\N^\infty)
=Q^{-1}\^m_\phi([\om]) \endaligned 
$$
and  
$$\aligned
\^m_\phi(\sg^{-n}([\om])) 
&= \sum_{\tau\in \N^n}\^m_\phi([\tau\om])
= \sum_{\tau\in \N^n}
\int\exp\(S_{\tau\om}(\phi-\P(\phi)|\tau\om|\)\,dm_\phi \\
&\le \sum_{\tau\in \N^n}||\exp\(S_\tau(\phi-\P(\phi)|\tau|\)||_0
\int\exp\(S_\om(\phi)-\P(\phi)|\om|\)\,dm_\phi \\
&=\exp\(S_\om(\phi)-\P(\phi)|\om|\)\,dm\sum_{\tau\in
\N^n}||\exp\(S_\tau(\phi)-\P(\phi)|\tau|\)||_0 \\
&\le Q\^m_\phi([\om]). \endaligned
$$
Let now $L$ be a Banach limit defined on the Banach space of 
all bounded sequences of real numbers. We define $\mu([\om]) 
=L\((\^m_\phi(\sg^{-n}([\om])))_{n\ge 0}\)$. Hence 
$Q^{-1}\^m_\phi([\om])\le \mu([\om])\le Q\^m_\phi([\om])$ and therefore
it is not difficult to check that the formula $\mu(A) 
=L\((\^m_\phi(\sg^{-n}(A)))_{n\ge 0}\)$ defines a finite non-zero
finitely additive measure on Borel sets of $\N^\infty$ satisfying 
$Q^{-1}\^m_\phi(A)\le \mu(A)\le Q\^m_\phi(A)$. Using now a theorem of
Calderon (Theorem 3.13 of [Fr]) and its proof one constructs a Borel
probability ($\sg$-additive) measure $\^\mu_\phi$ on $\N^\infty$
satisfying the formula 
$$
Q^{-1}\^m_\phi(A)\le \^\mu_\phi(A)\le Q\^m_\phi(A)
$$
for every Borel set $A\sbt\N^\infty$ with, perhaps, a larger constant $Q$.
Thus, to complete the proof 
of our theorem we only need to show total ergodicity of $\^\mu_\phi$ or
equivalently of $\^m_\phi$. Toward this end take a Borel set $A\in
\N^\infty$ with $\^m_\phi(A)>0$. Since 
the nested family of sets $\{[\tau]:\tau\in \N^*\}$ generates the
Borel $\sg$-algebra on $\N^\infty$, for every $n\ge 0$ and every
$\om\in \N^n$ we can find a 
subfamily $Z$ of $\N^*$ consisting of mutually incomparable words and
such that 
$A\sbt \bu\{[\tau]:\tau\in Z\}$ and $\sum_{\tau\in Z}\^m_\phi([\om\tau])\le 
2\^m_\phi(\om A)$, where $\om A=  
\{\om\rho:\rho\in A\}$. Then
$$\aligned
\^m_\phi\(\sg^{-n}(A)\cap [\om]\) 
&= \^m_\phi(\om A)\ge{1\over 2}\sum_{\tau\in Z}\^m_\phi([\om\tau]) \\
&= {1\over 2}\sum_{\tau\in Z} \int\exp
\(S_{\om\tau}(\phi-\P(\phi)|\om\tau|\)\,dm_\phi \\
&\ge {1\over 2Q}\exp\(S_\om(\phi-\P(\phi)|\om|\)||_0
    \sum_{\tau\in Z} \int\exp\(S_\tau(\phi-\P(\phi)|\tau|\)\,dm_\phi \\
&\ge {1\over
    2Q}\int\exp\(S_\om(\phi-\P(\phi)|\om|\)\,dm_\phi\sum_{\tau\in Z} 
\^m_\phi([\tau]) \\
&\ge {1\over 2Q}\^m_\phi([\om])\^m_\phi\((\bu\{[\tau]:\tau\in Z\}\)
\ge  {1\over 2Q}\^m_\phi(A)\^m_\phi([\om]).  \endaligned 
$$
Therefore 
$$
\aligned
\^m_\phi\(\sg^{-n}(\N^\infty\sms A)\cap [\om]\)
&=\^m_\phi\([\om]\sms \sg^{-
n}(A)\cap [\om]\) =\^m_\phi([\om])-\^m_\phi\(\sg^{-n}(A)\cap [\om]\) \\
&\le\(1-(2Q)^{-1}\^m_\phi(A)\)\^m_\phi([\om]).
\endaligned
$$ 
Hence for every Borel set $A\sbt \N^\infty$ with 
$\^m_\phi(A)<1$, for every $n\ge 0$, and for every $\om\in \N^n$ we get
\begin{equation}
\^m_\phi(\sg^{-n}(A)\cap[\om]\) \le \(1-(2Q)^{-1}(1-\^m_\phi(A))\)
\^m_\phi([\om]). 
\end{equation}
In order to conclude the proof of total ergodicity of $\sg$ suppose that 
$\sg^{-r}(A)=A$ for some integer $r\ge 1$ and some Borel set $A$ with
$0<\^m_\phi(A)<1$. Put $\g=1-(2Q)^{-1}(1-\^m_\phi(A))$. Note that $0<\g<1$.
In view of 
(2.5), for every $\om\in (\N^r)^*$ we get $\^m_\phi(A\cap[\om])=\^m_\phi(\sg^{-
|\om|}(A)\cap[\om]\) \le \g\^m_\phi([\om])$. Take now $\eta>1$ so small that 
$\g\eta<1$ and choose a subfamily $R$ of $(\N^r)^*$ consisting of mutually 
incomparable words and such that $A\sbt \bu\{[\om]:\om\in R\}$ and 
$\^m_\phi\(\bu\{[\om]:\om\in R\}\) \le \eta\^m_\phi(A)$. Then
$\^m_\phi(A)\le \sum_{\om\in 
R}\^m_\phi(A\cap[\om])\le \sum_{\om\in R}\g\^m_\phi([\om])=\g\^m_\phi
\(\bu\{[\om]:\om\in 
R\}\) \le \g\eta\^m_\phi(A)<\^m_\phi(A)$. This contradiction finishes the
proof. \end{proof}

\

\fr Set $\nu=\^m_\phi$. Clearly our assumption $\phi_i(X)\cap
\phi_j(X)=\es$ for $i,j\in\N$, $i\ne j$ implies that $\pi:\Sg\to J$ is
a homeomorphism; in particuluar, in view of Lemma~2.6, it establishes a
measure  preserving isomorphism between measure spaces $(\Sg,\nu)$ and 
$(J,m_\phi)$. To check that ${\it L}_\phi^*(\nu)=\l\nu$ take $g\in
C_b(\Sg)$ and compute
$$\aligned
\int g d{\it L}_0^*(\nu)
&=\int {\it L}_0(g)d\nu
 =\int {\it L}_0(g)(\pi^{-1}(x))d\nu\circ \pi^{-1}(x)
 =\int {\it L}_0(g)(\pi^{-1}(x))dm_\phi \\
&=\int\sum_{\tau\in\sg^{-1}(\pi^{-1}(x))}\exp(\phi(\tau)-\P(\phi))dm_\phi
    \\
&=\int\sum_{i\in\N}\exp\(\phi^{(i)}(x)-\P(\phi)\)g\circ\pi^{-1}(\phi_i(x)) 
  dm_\phi(x) \\
&=\int \^L_0(g\circ\pi^{-1})dm_\phi
 =\int g\circ\pi^{-1}dm_\phi
 =\int gd\nu. \endaligned
$$
Thus ${\it L}_0(\nu)=\nu$ and by the definition of ${\it L}_0$ and 
${\it L}_0^*$, ${\it L}_\phi^*(\nu)=\l\nu$. The fact that ${\it L}_\phi(h)=\l 
h$ follows
immediately from the definition of the operator ${\it L}_0$ and Theorem~2.7,
where $h=d\^\mu_\phi /d\^m_\phi$. Theorem~2.7 also implies that $h$ is 
bounded away from zero and infinity. In order to obtain H\"older continuity 
of the function $h$ and two convergence statements claimed in
Theorem~1.2 one may argue as follows: A well-known computation (see
[DU], comp [MU1]) shows that ${\it L}_0$ acts on the Banach space of bounded
uniformly continuous functions on $\N^\infty$ as an almost periodic operator 
(see [Ly], comp. [DU] and [MU1]). Using Theorem~2.7 and the theory of 
positive operators on lattices (see [Sc]) one then proves as in [DU]
that 1 is the only spectral point of modulus 1 and additionally that 1
is a simple eigenvalue of ${\it L}_0$. These facts and almost periodicity 
imply the first convergence statement of Theorem~1.2 and uniform
continuity of $h$. A similar computation as above produces constants
$0<\g<1$, $n\ge 1$ and $C\ge 0$ such that 
$$
||{\it L}_0^n(\g)||_\b\le C||\g||_0  + \g||g||_\b,
$$
where $||\g||_\b=V_\b(\g)+||g||_0$. This is so called the Ionescu-Tulcea
and Marinescu inequality. Using this inequality and Theorem~2.4 one checks
that the assumptions of the theorem of Ionescu-Tulcea and Marinescu (see
[IM], comp. [PU]) are satisfied. This theorem gives a nice spectral
decomposition of the operator ${\it L}_0$ acting on the space ${\it H}_\b$
of bounded H\"older continuous functions of order $\b$. Having this, a
relatively straightforward reasoning (comp. [PU]) shows H\"older
continuity of $h$ and the second convergence statement of Theorem~1.2.

\section{Equilibrium states} 
In this section we further investigate the
$\sg$-invariant measure $\^\mu_\phi$ introduced in Theorem~2.7. We
begin with the following technical result.
\begin{lemma} 
The following 3 conditions are equivalent
(a) $\int-\phi d\^\mu_\phi<\infty.$ 

(b)$\sum_{i\in\N}\inf(-\phi|_{[i]})\exp(\inf\phi|_{[i]})<\infty.$

(c) $\H_{\^\mu_\phi}(\a)<\infty$, where $\a=\{[i]:i\in \N\}$ is the
partition of $\Sg$ into initial cylinders of length 1.
\end{lemma}
\begin{proof} Suppose that $\int-\phi d\^\mu_\phi<\infty.$ It
means that $\sum_{i\in \N}\int_{[i]}-\phi d\^\mu_\phi<\infty$ and
consequently 
$$\aligned
\infty
&>\sum_{i\in\N}\inf(-\phi|_{[i]}) \int_{[i]} d\^\mu_\phi 
 =\sum_{i\in\N}\inf(-\phi|_{[i]}) \int_{[i]} hd\^m_\phi \\
&\ge Q^{-1}\sum_{i\in\N}\inf(-\phi|_{[i]}) \^m_\phi([i])
 =Q^{-1}\sum_{i\in\N}\inf(-\phi|_{[i]})
  \int_X\exp(\phi^{(i)}(x)-\P(\phi))dm_\phi(x) \\
&=Q^{-1}\ep^{-\P(\phi)}\sum_{i\in\N}\inf(-\phi|_{[i]})\int_X
  \exp(\phi^{(i)}(x))dm_\phi(x) \endaligned
$$
Thus 
$$\aligned
\infty
&>\sum_{i\in\N}\inf(-\phi|_{[i]})\int_X\exp(\phi^{(i)}(x))dm_\phi(x) 
\ge \sum_{i\in\N}\inf(-\phi|_{[i]})\exp(\inf_X(\phi^{(i)}) \\
&=\sum_{i\in\N}\inf(-\phi|_{[i]})\exp(\inf\phi|_{[i]})\endaligned
$$
Now suppose that $\sum_{i\in\N}\inf(-\phi|_{[i]}) 
\exp(\inf\phi|_{[i]})<\infty$. We shall show that
$\H_{\^\mu_\phi}(\a)<\infty$. So,
$$
\H_{\^\mu_\phi}(\a)
=\sum_{i\in\N}-\^\mu_\phi([i])\log\^\mu_\phi([i])
\le \sum_{i\in\N}-Q\^m_\phi([i])\(\log\^m_\phi([i]) -\log Q\).
$$
But $\sum_{i\in\N}-Q\^m_\phi([i])( -\log Q)=Q\log Q$, so it
suffices to show that
$$
\sum_{i\in\N}-\^m_\phi([i])\log\^m_\phi([i]) <\infty.
$$
But 
$$\aligned
\sum_{i\in\N}-\^m_\phi([i])\log\^m_\phi([i])
&=\sum_{i\in\N}-\^m_\phi([i])\log\lt(\int_X\exp\(\phi^{(i)}-
\P(\phi)\)\rt)dm_\phi \\
&\le \sum_{i\in\N}-\^m_\phi([i])(\inf_X\phi^{(i)}-\P(\phi)). \endaligned
$$
But $\sum_{i\in\N}\^m_\phi([i])\P(\phi)=\P(\phi)$, so it suffices to
show that $\sum_{i\in\N}-\^m_\phi([i])\inf_X\phi^{(i)} <\infty$. And
indeed, using Lemma~2.1 we get
$$
\sum_{i\in\N}-\^m_\phi([i])\inf_X\phi^{(i)} 
= \sum_{i\in\N}\^m_\phi([i])\sup_X(-\phi^{(i)})
\le \sum_{i\in\N}\^m_\phi([i])\(\inf_X(-\phi^{(i)})+\log Q\).
$$
Since $\sum_{i\in\N}\^m_\phi([i])\log Q=\log Q$, it is enough to show
that 
$$
\sum_{i\in\N}\^m_\phi([i])\inf_X(-\phi^{(i)}) <\infty.
$$ 
And indeed, 
$$
 \sum_{i\in\N}\^m_\phi([i])\inf_X(-\phi^{(i)})
=\sum_{i\in\N}\int\exp(\phi^{(i)}-\P(\phi))dm_\phi\inf_X(-\phi^{(i)}) 
$$
But in view of (1.4) $\phi^{(i)}$ are negative everywhere for all $i$
large enough, say $i\ge k$. Then using Lemma~2.1 again we get
$$
 \sum_{i\ge k}\^m_\phi([i])\inf_X(-\phi^{(i)})
\le
 \ep^{-\P(\phi)}Q\sum_{i\ge k}\exp(\inf_X(\phi^{(i)}))\inf_X(-\phi^{(i)}) 
$$
which is finite due to our assumption. Hence,
$\sum_{i\in\N}\^m_\phi([i])
\inf_X(-\phi^{(i)})<\infty$. Finally suppose that
$\H_{\^\mu_\phi}(\a)<\infty$. We need to show that $\int-\phi
d\^\mu_\phi <\infty.$ We have
$$
\infty
>\H_{\^\mu_\phi}(\a)
= \sum_{i\in\N}-\^m_\phi([i])\log\(\^m_\phi([i])\)
\le \sum_{i\in\N}-\^m_\phi([i])\(\inf(\phi|_{[i]}-\P(\phi)-\log Q)\).
$$
Hence $\sum_{i\in\N}-\^m_\phi([i])\inf(\phi|_{[i]})<\infty$ and
therefore 
$$
\int-\phi d\^\mu_\phi
= \sum_{i\in\N}\int_{[i]} -\phi d\^\mu_\phi
\le \sum_{i\in\N}\sup(-\phi|_{[i]})\^m_\phi([i])
= \sum_{i\in\N}-\inf(\phi|_{[i]})\^m_\phi([i])
<\infty.
$$
The proof is complete. \end{proof}

\

\fr By Theorem 3 of [Sa] we know that $\sup\{\hmu(\sg)+\int\phi
d\mu\}=\P(\phi)$, where the supremum is taken over all $\sg$-invariant
probability measures such that $\int-\phi d\mu<\infty$. We call a
$\sg$-invariant probability measure $\mu$ an equilibrium state of the
potential $\phi$ if $\hmu(\sg)+\int\phi d\mu=\P(\phi)$. We shall prove
the following.

\begin{theorem} 
If $\sum_{i\in\N}\inf(-\phi|_{[i]})\exp(\inf\phi|_{[i]})<\infty$, then
$\^\mu_\phi$ is an equilibrium state of the potential $\phi$
satisfying $\int-\phi d\^\mu_\phi<\infty$. 
\end{theorem}
\begin{proof} It follows from Lemma~3.1 that $\int-\phi
d\^\mu_\phi <\infty$. To show that $\^\mu_\phi$ is an equilibrium
state of the potential $\phi$ consider $\a=\{[i]:i\in\N\}$, the
partition of $\Sg$ into initial cylinders of length one. By Lemma~3.1,
$\H_{\^\mu_\phi}(\a)<\infty$. Applying the Breiman-Shanon-McMillan
theorem and the Birkhoff ergodic theorem we therefore get for
$\^\mu_\phi$-a.e. $\om\in\Sg$ 
$$\aligned
\text{h}_{\mu_\phi}(\sg)
&\ge\text{h}_{\mu_\phi} (\sg,\a)
=\lim_{n\to\infty}{-1\over n} \log([\om|_n]) \\
&=\lim_{n\to\infty}{-1\over n}
\log\lt(\int\exp\(S_\om(\phi)(x)d\mu_\phi -\P(\phi)n\)\rt) \\
&=\lim_{n\to\infty}{-1\over n} \log\lt(\int\exp\(\sum_{j=0}^{n-1}\phi 
  (\sg^j(\om|_n\tau))d\mu_\phi(\tau) -\P(\phi)n\)\rt) \\
&\ge \limsup_{n\to\infty}{-1\over n} \log\lt(\int\exp\(\sum_{j=0}^{n-1}\phi 
  (\sg^j(\om))+\log Q -\P(\phi)n\)\rt) \\
&=\lim_{n\to\infty}{-1\over n} \sum_{j=0}^{n-1}\phi(\sg^j(\om))
+\P(\phi)
=-\int\phi d\^\mu_\phi+\P(\phi). \endaligned
$$
Hence $\text{h}_{\mu_\phi}(\sg)+\int\phi d\^\mu_\phi\ge \P(\phi)$,
which in view
of the variational principle (see Theorem 3 in [Sa]), implies that
$\^\mu_\phi$ is an equilibrium state for the potential $\phi$. The
proof is finished. \end{proof}

\

\fr{\bf Acknowledgment.} We would like to thank Prof. Sarig Omri who has 
read the preliminary version of this paper for his helpful
remarks and comments which improved the present version. In particular we 
would like to thank him for a direct proof of Lemma~1.1.

\

\bibliographystyle{amsalpha}

\begin{thebibliography}{A}


\bibitem [Bo]{}
R. Bowen, Equilibrium States and Ergodic Theory of Anosov
Diffeomrphisms, L. N. Math. 470, Berlin, Heidelberg, New York, 
Springer-Verlag, (1975).

\bibitem [DU]{}
M. Denker, M. Urba\'nski, Ergodic theory of equilibrium states
for rational maps, Nonlinearity 4, (1991), 103-134.

\bibitem [Fr]{}
N. Friedman, Introduction to Ergodic Theory, New York, Cincinati,
Toronto, London, Melbourne, Van Nostrand Reinhold Company, (1970).

\bibitem [IM]{}
C. Ionescu-Tulcea, G. Marinescu, Th\'eorie regodique pour des classes
d'operations non-complement continues, Ann. Math. 52, (1950), 140-147.

\bibitem [Ly]{}
M. Lyubich, Entropy properties of rational endomorphisms of the Riemann
sphere, Ergod. Th. and Dynam. Sys. 3, (1983), 351-386.

\bibitem [MU1]{}
R.D. Mauldin, M. Urba\'nski, Dimensions and measures in infinite
iterated 
function systems, Proc. London Math. Soc. (3) 73(1996), 105-154.     
                
\fr [MU2]{}
R.D. Mauldin, M. Urba\'nski, Parabolic iterated function systems, Preprint 
1998.   

\bibitem [PU]{}
F. Przytycki, M. Urba\'nski, Fractals in the Plane - Ergodic Theory
Methods, to appear.

\bibitem[Sa]{}
O. M. Sarig, Theormodynamic formalism for countable Markov shifts,
to appear Ergod. Th. and Dynam. Sys..

\bibitem [Sc]{}
H. Schaefer, Banach Lattices and Positive Operators. (1974), Springer.
 
\bibitem [Wa]{} 
P. Walters, Invariant measures and equilibrium states for some mappings which 
expand distances, Transactions of A.M.S. 236 (1978), 121 - 153.

\end{thebibliography}

\end{document}